\documentclass[10pt,a4paper,twoside,final]{amsart}

\usepackage[latin1]{inputenc}
\usepackage[T1]{fontenc}
\usepackage{amssymb}
\usepackage[UKenglish]{babel}

\numberwithin{equation}{section}

\newtheorem{theorem}{Theorem}
\newtheorem{lemma}[theorem]{Lemma}
\newtheorem{prop}[theorem]{Proposition}
\newtheorem{cor}[theorem]{Corollary}
\newtheorem{claim}{Claim}

\theoremstyle{remark}
\newtheorem{remark}[theorem]{Remark}

\renewcommand{\Re}{\operatorname{Re}}

\newcommand{\abs}[1]{\lvert#1\rvert}
\newcommand{\bigabs}[1]{\bigl\lvert#1\bigr\rvert}
\newcommand{\biggabs}[1]{\biggl\lvert#1\biggr\rvert}
\newcommand{\norm}[1]{\lVert#1\rVert}
\newcommand{\bignorm}[1]{\bigl\lVert#1\bigr\rVert}
\newcommand{\biggnorm}[1]{\biggl\lVert#1\biggr\rVert}

\newcommand{\BMOA}{\mathit{BMOA}}
\newcommand{\VMOA}{\mathit{VMOA}}
\newcommand{\VMO}{\mathit{VMO}}
\newcommand{\BMO}{\mathit{BMO}}
\newcommand{\Bloch}{\mathcal{B}}

\newcommand{\D}{\mathbb{D}}
\newcommand{\T}{\mathbb{T}}
\newcommand{\C}{\mathbb{C}}

\title{Compact and weakly compact composition operators on BMOA}

\author[Laitila]{Jussi Laitila}
\address{Institute for Social and Economic Research, University of
Essex, Colchester CO4 3SQ, United Kingdom}
\email{jlaitila@essex.ac.uk}

\author[Nieminen]{Pekka J.\ Nieminen}
\address{Department of Mathematics and Statistics, University of
Helsinki, PO~Box~68, FI-00014 Helsinki, Finland}
\email{pjniemin@cc.helsinki.fi}

\author[Saksman]{Eero Saksman}
\address{Department of Mathematics and Statistics, University of
Helsinki, PO~Box~68, FI-00014 Helsinki, Finland}
\email{eero.saksman@helsinki.fi}

\author[Tylli]{Hans-Olav Tylli}
\address{Department of Mathematics and Statistics, University of
Helsinki, PO~Box~68, FI-00014 Helsinki, Finland}
\email{hojtylli@cc.helsinki.fi}

\thanks{The first author was supported by the Academy of
Finland, project 118422. The second and third authors were
supported by the Finnish CoE in Analysis and Dynamics Research,
and by the Academy of Finland, projects 118422 \& 126420 and
113826 \& 118765, respectively.}


\date{20 January 2010}

\begin{document}

\begin{abstract}
Any analytic map $\varphi$ of the unit disc $\D$ into itself induces
a composition operator $C_\varphi$ on $\BMOA$, mapping
$f \mapsto f \circ \varphi$, where $\BMOA$ is the Banach space
of analytic functions $f\colon \D \to \C$ whose boundary values have
bounded mean oscillation on the unit circle. We show that
$C_\varphi$ is weakly compact on $\BMOA$ precisely when it is compact on
$\BMOA$, thus solving a question initially posed by Tjani and by Bourdon,
Cima and Matheson in the special case of $\VMOA$.
As a crucial step of our argument we simplify the compactness criterion
due to Smith for $C_\varphi$ on $\BMOA$ and show that 
his condition on the Nevanlinna counting function alone characterizes
compactness. Additional equivalent compactness criteria are established.
Furthermore, we prove the unexpected result that compactness of 
$C_\varphi$ on $\VMOA$ implies compactness even from the Bloch space
into $\VMOA$.
\end{abstract}

\maketitle

\section{Introduction}

Let $\D$ be the open unit disc of the complex plane $\C$.
The space $\BMOA$ consists of the analytic functions
$f\colon \D \to \C$ whose boundary values have \emph{bounded
mean oscillation} on the unit circle $\T$. Equivalently, $f$ belongs to
$\BMOA$ if and only if the seminorm
\[
   \abs{f}_* = \sup_{a\in\D} \bignorm{f\circ\sigma_a-f(a)}_{H^2} 
\]
is finite, where $\norm{\cdot}_{H^2}$ is the standard norm of the Hardy
space $H^2$ and $\sigma_a(z) = (a-z)/(1-\bar{a}z)$ is the automorphism
of $\D$ that exchanges the points $0$ and $a$. Then $\BMOA$
becomes a Banach space under the norm
$\norm{f}_* = \abs{f(0)} + \abs{f}_*$. Furthermore, $\VMOA$ is the
closed subspace of $\BMOA$ consisting of those functions $f$ whose
boundary values have \emph{vanishing mean oscillation}, or equivalently,
which satisfy
\[
   \lim_{\abs{a}\to 1} \bignorm{f\circ\sigma_a-f(a)}_{H^2} = 0.
\]
We refer e.g.\ to \cite{Girela} and \cite{Zhu} for more
information on the spaces $\BMOA$ and $\VMOA$.

If $\varphi\colon \D \to \D$ is an analytic map,
then the \emph{composition operator} $C_\varphi$ induced by $\varphi$ is
the linear map defined by $C_\varphi f = f\circ\varphi$ for all
analytic functions $f\colon \D \to \C$. It is well known that
$C_\varphi$ is always bounded from $\BMOA$ into itself and that
$C_\varphi$ preserves $\VMOA$ if and only if $\varphi \in \VMOA$; see
e.g.\ \cite{Stephenson}, \cite{AFP} and \cite{BCM}.
Composition operators have been intensively studied on various spaces 
of analytic functions, and we refer to \cite{CowenMacCluer} or
\cite{Shapiro} for more about the classical background.

Recall that a linear operator is \emph{compact} if it takes bounded
sets into sets having a compact closure. The compactness of a
composition operator $C_\varphi$ acting on $\BMOA$ (or on its subspace
$\VMOA$) has been investigated by several authors and various kinds of
characterizations are known; see e.g.\ \cite{Tjani}, \cite{BCM},
\cite{Smith}, \cite{MakhTjani}, \cite{WirthsXiao},
\cite{Laitila:vector}, \cite{Wulan}, \cite{Laitila:weighted} and
\cite{WZZ}. In particular, Smith \cite{Smith} proved that $C_\varphi$ 
is compact on $\BMOA$ if and only if $\varphi$ satisfies
the following pair of conditions:
\begin{equation}\tag{S1}\label{Sm1}
   \lim_{\abs{\varphi(a)} \to 1} \sup_{0 < \abs{w} < 1} \abs{w}^2 
   N(\sigma_{\varphi(a)} \circ \varphi \circ \sigma_a,w) = 0,
\end{equation}
and for all $0 < R < 1$,
\begin{equation}\tag{S2}\label{Sm2}
   \lim_{t\to 1} \sup_{\{a: \abs{\varphi(a)} \leq R\}} 
   \bigabs{ \{\zeta \in \T:
   \abs{(\varphi \circ \sigma_a)(\zeta)} > t\} } = 0.
\end{equation}
Above $N(\psi,\cdot)$ denotes the Nevanlinna counting function of an
analytic self-map $\psi$ of the disc, $\varphi(\zeta)$ is the
radial limit of $\varphi$ for a.e.\ $\zeta$ on the unit circle $\T$,
and $\abs{E}$ stands for the normalized Lebesgue measure of
sets $E \subset\T$.
Recently the first author \cite{Laitila:weighted}
showed that \eqref{Sm1} is equivalent to the condition
\begin{equation}\tag{L}\label{L}
   \lim_{\abs{\varphi(a)} \to 1}
   \bignorm{\sigma_{\varphi(a)}\circ\varphi\circ\sigma_a}_{H^2} = 0,
\end{equation}
which is technically more convenient for our later purposes.

A well-known open problem concerning composition operators
is that of characterizing the \emph{weak compactness} of $C_\varphi$
on the non-reflexive spaces $\BMOA$ and $\VMOA$. Recall that an
operator is weakly compact provided it takes bounded sets into sets
whose closure is compact in the weak topology of the space. For
$C_\varphi$ acting on $\VMOA$ this problem was explicitly
posed in \cite{Tjani} and \cite{BCM}, and for the $\BMOA$ case
it was stated in \cite{Laitila:wcompact,Laitila:vector}.
Partial results for $\VMOA$ were obtained in \cite{MakhTjani} and
\cite{CimaMatheson}. For instance, if $\varphi \in \VMOA$ and
$\varphi(\D)$ is contained in a polygon inscribed in $\overline{\D}$
\cite[Cor.~5.4]{MakhTjani}, or if $\varphi$ is univalent
\cite[p.~940]{CimaMatheson}, then compactness and weak compactness are
equivalent for $C_\varphi$ on $\VMOA$. It is natural to
conjecture that the same equivalence should persist for arbitrary
symbols $\varphi$ even on $\BMOA$, especially because
a similar phenomenon is known to occur for composition operators
on many other classical non-reflexive spaces, such as
$H^1$ \cite{Sarason}, $H^\infty$ (see e.g.\ \cite{AGL}) and
Bloch spaces \cite{MadiganMatheson,LST}.

In the present paper we provide a solution to the above problem.
Our main result reads as follows:

\begin{theorem} \label{thm:Main1}
Let $\varphi\colon \D \to \D$ be an analytic map.
Then the following conditions are equivalent:
\begin{enumerate}
\renewcommand{\theenumi}{\roman{enumi}}
\item
$C_\varphi\colon \BMOA \to \BMOA$ is compact.
\item\label{it:WCpt}
$C_\varphi\colon \BMOA \to \BMOA$ is weakly compact.
\item 
$\varphi$ satisfies condition \eqref{Sm1}.
\item\label{it:L}
$\varphi$ satisfies condition \eqref{L}.
\end{enumerate}
\end{theorem}

A key ingredient of our argument is the surprising result that
condition \eqref{L} (and consequently also \eqref{Sm1}) actually implies
\eqref{Sm2}. This result is proved in Section~\ref{sec:Cpt}.
Thus our work substantially clarifies and simplifies the existing
compactness criteria for composition operators on $\BMOA$. 
The proof of Theorem~\ref{thm:Main1} is then completed by verifying that
\eqref{it:WCpt} implies~\eqref{it:L}. This step is carried out in
Section~\ref{sec:Weak}, where the argument is based on an idea of
Le{\u\i}bov \cite{Leibov} on how to construct explicit isomorphic
copies of the sequence space $c_0$ inside $\VMOA$.

As a by-product the results of Section~\ref{sec:Cpt} answer a
recent question of Wulan, Zheng and Zhu \cite{WZZ}. Namely, it follows
that the condition 
$\lim_{\abs{a} \to 1} \abs{\sigma_a \circ \varphi}_* = 0$
is sufficient for the compactness of $C_\varphi$ on $\BMOA$.
The necessity was earlier observed by Wulan~\cite{Wulan}.

In Section~\ref{sec:MeanOsc} we further reformulate \eqref{L} as a 
pseudo-hyperbolic mean oscillation condition for the 
boundary values of the symbol as follows:
\begin{equation}\tag{A}\label{aver}
   \frac{1}{\abs{I}^2} \int_I\int_I
   \rho\bigl( \varphi(\zeta),\varphi(\xi) \bigr)^2
   \,\abs{d\zeta}\abs{d\xi} \to 0
   \quad\text{as}\quad
   \biggabs{ \frac{1}{\abs{I}}
   \int_I\varphi(\zeta) \,\abs{d\zeta} } \to 1.
\end{equation}
Here $\rho$ denotes the pseudo-hyperbolic metric, $I \subset \T$ is a
boundary arc and the integration is with respect to the normalized
Lebesgue measure on $\T$. 

Section~\ref{sec:VMOA} collects together some related results in the
$\VMOA$ setting. We observe that the analogue of
Theorem~\ref{thm:Main1} holds on $\VMOA$ (that is, for symbols
$\varphi \in \VMOA$), where \eqref{L} can be replaced by
$\lim_{\abs{a} \to 1}
\bignorm{\sigma_{\varphi(a)}\circ\varphi\circ\sigma_a}_{H^2} = 0$.
Moreover, we prove that one may substitute the genuine hyperbolic
metric for the pseudo-hyperbolic metric in the $\VMOA$ version of
condition \eqref{aver}. As an unexpected consequence, $C_\varphi$
is compact on $\VMOA$ if and only if it is compact from the Bloch
space to $\VMOA$.

\section{Compactness characterization}
\label{sec:Cpt}

In this section we prove that condition~\eqref{L} alone is enough
to characterize the compactness of $C_\varphi$ on the space $\BMOA$.
It is known that \eqref{L} is equivalent to Smith's first
condition~\eqref{Sm1}; see \cite[Remark~3.3]{Laitila:weighted}. Thus,
in view of Smith's compactness criterion
consisting of the pair \eqref{Sm1} and \eqref{Sm2}, our work reduces to
showing that \eqref{Sm2} is actually implied by~\eqref{Sm1}, or
by~\eqref{L}:

\begin{theorem} \label{thm:Superfluous}
Condition \eqref{L} implies \eqref{Sm2} for any analytic map
$\varphi\colon \D \to \D$. Hence $C_\varphi\colon \BMOA \to \BMOA$ is
compact if and only if \eqref{L} holds.
\end{theorem}

We mostly work with \eqref{L} because it is technically very convenient
for our arguments and also allows for quite appealing
reformulations in terms of the boundary values of $\varphi$.
In particular, by expressing the $H^2$ norm as an $L^2$ norm on $\T$ and
performing a change of variable using the automorphism $\sigma_a$, we
get
\begin{equation} \label{eq:POsc} \begin{split}
   \bignorm{\sigma_{\varphi(a)}\circ\varphi\circ\sigma_a}_{H^2}^2
   &= \int_\T \rho\bigl( \varphi(\sigma_a(\zeta)),\varphi(a) \bigr)^2
      \,\abs{d\zeta}  \\
   &= \int_\T \rho\bigl( \varphi(\zeta), \varphi(a) \bigr)^2
      P_a(\zeta) \,\abs{d\zeta},
\end{split} \end{equation}
where $P_a(\zeta) = (1-\abs{a}^2)/\abs{\zeta-a}^2$ is the Poisson
kernel for $a \in \D$ and
$\rho(z,w) = \abs{z-w}/\abs{1-\overline{w}z}$ denotes the
pseudo-hyperbolic distance in $\overline{\D}$ (observe that $\rho$
extends to the boundary $\T$ in a natural way if we agree that
$\rho(z,z) = 0$ for $z\in\T$).
Thus \eqref{L} can be seen as a kind of vanishing mean
oscillation condition with respect to the pseudo-hyperbolic metric. We
will elaborate on this point further in Section~\ref{sec:MeanOsc}.

It is useful to observe that if $\varphi$ 
satisfies condition \eqref{L}, or equivalently \eqref{Sm1},
then one has $\abs{\varphi} < 1$ a.e.\ on $\T$. This can be checked 
by a straightforward density point argument. 

The proof of Theorem~\ref{thm:Superfluous} depends on the following
lemma, which exhibits a uniform density estimate for Lebesgue measurable
sets on $\T$. Since we have been unable to find a reference for this
kind of result, we include a proof. Here and elsewhere in the text we
use the following notation for closed arcs of $\T$:\ when
$re^{i\theta} \in \D$ with $0 \leq r < 1$, set
\[
   I(re^{i\theta}) = \{ e^{it} : \abs{t-\theta} \leq \pi(1-r) \}.
\]
Thus $I(re^{i\theta})$ denotes the arc of $\T$ whose midpoint is
$e^{i\theta}$ and (normalized) length $\abs{I(re^{i\theta})} = 1-r$.

\begin{lemma} \label{le:Density}
Suppose that $E \subset \T$ is a measurable set with $\abs{E} > 0$.
Then there is a measurable set $E' \subset E$ such that $\abs{E'} > 0$
and
\[
   \frac{\abs{I(r\zeta) \cap E}}{\abs{I(r\zeta)}}
   \geq \frac{1}{8}\abs{E}
\]
for every $0 \leq r < 1$ and $\zeta \in E'$.
\end{lemma}

\begin{proof}
We say that an arc $I(a) \subset \T$ is \emph{dyadic} if
$a = (1-2^{-n})\exp(2\pi ik/2^n)$ for some $n \geq 0$ and
$0 \leq k \leq 2^n-1$. Note that any pair of dyadic arcs either
are nested or have disjoint interiors.

Put $\lambda = 1-\tfrac{1}{2}\abs{E} \in
\mathopen[\tfrac{1}{2},1\mathclose)$, and let
$C$ be the set of all points $\zeta \in \T$ for which there
exists a dyadic arc $I$ containing $\zeta$ and satisfying
$\abs{I\cap E^c} > \lambda \abs{I}$, where $E^c = \T\setminus E$.
Since for each $\zeta \in C$ there is a maximal one (in terms
of inclusion) among such dyadic arcs, we have
$C = \bigcup_j I_j$, where the $I_j$'s are
dyadic arcs with disjoint interiors satisfying
$\abs{I_j\cap E^c} > \lambda \abs{I_j}$. Summing over $j$ and
noting that almost every point of $E^c$ belongs to $C$ by the
Lebesgue density theorem, we then get
$\abs{E^c} = \abs{C \cap E^c} > \lambda\abs{C}$.
Thus $\abs{C} < \abs{E^c}/\lambda = (1-\abs{E})/\lambda < 1$,
and so $\abs{C^c} > 0$.

To finish the proof, note that for almost every
$\zeta \in C^c$ we have $\zeta \in E$ and also
$\abs{I \cap E} \geq (1-\lambda)\abs{I} = \tfrac{1}{2} \abs{E}\abs{I}$
for all dyadic arcs $I$ containing $\zeta$. Moreover, for every
arc $I(r\zeta) \subset \T$ there exists a dyadic arc $I$ such that
$\zeta \in I \subset I(r\zeta)$ and
$\abs{I} > \tfrac{1}{4}\abs{I(r\zeta)}$. These observations
prove the lemma with $E' = C^c \cap E$.
\end{proof}

\begin{proof}[Proof of Theorem~\ref{thm:Superfluous}]
As a preparatory step we first establish a M\"{o}bius-invariant
version of condition \eqref{L}. Let
$\varphi_b = \varphi\circ\sigma_b$ for $b \in \D$.
Then the following identity can be verified just by inspection
and using the self-inverse property of the automorphisms:
\[
   \sigma_{\varphi_b(a)} \circ \varphi_b \circ \sigma_a
   = \bigl[ \sigma_{\varphi(\sigma_b(a))} \circ \varphi \circ
            \sigma_{\sigma_b(a)} \bigr] \circ
     \bigl[ \sigma_{\sigma_b(a)} \circ \sigma_b \circ \sigma_a \bigr].
\]
Note that the composite mapping enclosed in the last brackets
is a disc automorphism that fixes the origin, hence a
rotation. Therefore
\[
   \bignorm{\sigma_{\varphi_b(a)} \circ \varphi_b \circ \sigma_a}_{H^2}
   = \bignorm{\sigma_{\varphi(\sigma_b(a))} \circ \varphi \circ
              \sigma_{\sigma_b(a)}}_{H^2}.
\]
Now, in view of \eqref{eq:POsc} and the fact that
$P_a(\zeta) \geq \frac{1}{4}\abs{I(a)}^{-1}$ for $\zeta \in I(a)$,
condition~\eqref{L} implies the following: Given $\varepsilon > 0$,
there exists $\eta < 1$ such that
\begin{equation} \label{eq:Invariant}
   \frac{1}{\abs{I(a)}}
   \int_{I(a)} \rho\bigl(\varphi_b(\zeta),\varphi_b(a)\bigr)^2
   \,\abs{d\zeta}
   \leq \varepsilon
\end{equation}
whenever $a$ and $b$ satisfy $\abs{\varphi_b(a)} \geq \eta$.

For the actual proof of Theorem~\ref{thm:Superfluous} we argue by
contradiction, assuming that \eqref{L} holds but \eqref{Sm2} does not.
Since \eqref{Sm2} fails, there are constants $R < 1$ and $c > 0$,
points $b_k \in \D$, and numbers $0 < t_k < 1$ with $t_k \to 1$ such
that for all $k \geq 1$ we have $\abs{\varphi(b_k)} \leq R$ and the
sets
\[
   E_k = \bigl\{\zeta \in \T: \text{the radial limit $\varphi_k(\zeta)$
      exists and $\abs{\varphi_k(\zeta)} > t_k$}\bigr\}
\]
satisfy
$\abs{E_k} \geq c$, where $\varphi_k = \varphi\circ\sigma_{b_k}$.
By Lemma~\ref{le:Density} we can further
find sets $E_k' \subset E_k$ such that $\abs{E_k'} > 0$ and
\begin{equation} \label{eq:Density}
   \frac{\abs{I(r\zeta) \cap E_k}}{\abs{I(r\zeta)}} \geq \frac{c}{8}
   \quad\text{for $0 \leq r < 1$, $\zeta \in E_k'$.}
\end{equation}

Let $\varepsilon = c/16$. We may choose $\eta$ large enough so that
$R < \eta < 1$ and \eqref{eq:Invariant} holds for
$\abs{\varphi_b(a)} \geq \eta$. Fix $k$ such that $t_k \geq \eta$.
Recall that by the definition of $E_k$ we have
$\abs{\varphi_k(r\zeta)} \to \abs{\varphi_k(\zeta)} > t_k$
as $r \to 1$ for each $\zeta \in E_k$. In particular, we
can fix a point $\zeta_k \in E_k'$ with this property. Moreover, since
$\abs{\varphi_k(0)} = \abs{\varphi(b_k)} \leq R$, it follows
from continuity that there is a radius $0 < r_k < 1$ such that
$\abs{\varphi_k(r_k\zeta_k)} = \eta$. Let $a_k = r_k\zeta_k$.
By elementary geometry it holds  for each $\zeta \in E_k$ that
$\rho(\varphi_k(\zeta),\varphi_k(a_k)) \geq
\rho(t_k,\eta)$. Hence we can use \eqref{eq:Density} to obtain the
estimate
\[
   \frac{1}{\abs{I(a_k)}} \int_{I(a_k)}
      \rho\bigl(\varphi_k(\zeta),\varphi_k(a_k)\bigr)^2 \,\abs{d\zeta}
   \geq \frac{\abs{I(a_k) \cap E_k}}{\abs{I(a_k)}} \rho(t_k,\eta)^2
   \geq \frac{c}{8} \rho(t_k,\eta)^2.
\]
Since this estimate holds for all sufficiently large $k$,
we may let $k \to \infty$. In this case $\rho(t_k,\eta) \to 1$, which
leads to a contradiction with \eqref{eq:Invariant} by the choice of
$\varepsilon$.
\end{proof}

\begin{remark}
It is appropriate to note that condition~\eqref{Sm2} alone
does not ensure the compactness of $C_\varphi$ on $\BMOA$.
For example, if $\varphi(z) = \tfrac{1}{2}(1+z)$, then one may
check that \eqref{Sm2} holds but $C_\varphi$ fails to be compact.
For instance, it is not difficult to see that
$\bignorm{\sigma_{\varphi(a)}\circ\varphi\circ\sigma_a}_{H^2} \to 1$
as $a \to 1$. We leave the details to the reader.
\end{remark}

We close this section by addressing a question recently posed by
Wulan, Zheng and Zhu~\cite{WZZ}. Based on an
earlier work by Wulan~\cite{Wulan}, they showed that the single
condition  
\begin{equation}\tag{W1}\label{W1}
   \lim_{n\to \infty} \abs{\varphi^n}_* = 0
\end{equation} 
characterizes the compactness of $C_\varphi$ on $\BMOA$.
The earlier result in \cite{Wulan} involved the additional
condition
\begin{equation}\tag{W2}\label{W2}
   \lim_{\abs{a} \to 1} \abs{\sigma_a \circ \varphi}_* = 0,
\end{equation}
and so it was asked in \cite[Sec.~4]{WZZ} whether
\eqref{W2} alone would suffice to characterize when $C_\varphi$ is
compact on $\BMOA$. This is indeed the case.

\begin{cor}
Let $\varphi\colon \D \to \D$ be an analytic map. Then $C_\varphi$ is
compact on $\BMOA$ if and only if \eqref{W2} holds.
\end{cor}

\begin{proof}
It is enough to observe that $\abs{\sigma_{\varphi(a)} \circ \varphi}_*
\geq \bignorm{\sigma_{\varphi(a)}\circ\varphi\circ\sigma_a}_{H^2}$,
whence \eqref{W2} trivially implies~\eqref{L}.
\end{proof}

\section{Weak compactness}
\label{sec:Weak}

After the work of the preceding section the only step that remains
to be proved in Theorem~\ref{thm:Main1} is that \eqref{it:WCpt}
implies~\eqref{it:L}.
Equivalently, if the map $\varphi$ fails to satisfy condition \eqref{L},
then we must show that the composition operator $C_\varphi$ is not
weakly compact on $\BMOA$. This will be accomplished separately in
Proposition~\ref{prop:Comp} below.

Our argument depends on the following proposition which is
essentially due to Le{\u\i}bov \cite{Leibov} and provides
information about the subspace structure of $\VMOA$ (cf.\
Remark~\ref{re:Sub} below). As usual, here $c_0$ denotes the
Banach space of complex sequences converging to zero endowed
with the supremum norm $\norm{\cdot}_\infty$. The proof given
below is an adaptation of Le{\u\i}bov's argument; he worked in the
space $\VMO(\T)$ on the unit circle, but we work directly in the disc.

\begin{prop} \label{prop:c0}
Let $(f_n)$ be a sequence in $\VMOA$ such that $\norm{f_n}_* = 1$ for
all $n$ and $\norm{f_n}_{H^2} \to 0$ as $n \to \infty$.
Then there exists a subsequence $(f_{n_k})$ which is equivalent
to the natural basis of $c_0$; that is, for which the map
$(\lambda_k) \mapsto \sum_k \lambda_k f_{n_k}$ is
an isomorphism from $c_0$ into $\VMOA$.
\end{prop}

\begin{proof}
For brevity we write
\[
   \gamma(f,a) = \bignorm{f\circ\sigma_a-f(a)}_{H^2}
\]
whenever $f \in H^2$ and $a \in \D$. Note that $\gamma(f,a)$ defines
a seminorm with respect to $f$ for each $a$. We also have
$\gamma(f,a) \leq \norm{f\circ\sigma_a}_{H^2} \leq c_a\norm{f}_{H^2}$
for some $c_a > 0$, where $c_a$ is an increasing function of $\abs{a}$.
Therefore
\[
   \sup \{\gamma(f_n,a) : \abs{a} \leq r\} \to 0
   \quad\text{as $n \to \infty$}
\]
for any $0 < r < 1$. On the other hand, the $\VMOA$ condition
says that $\gamma(f_n,a) \to 0$ as $\abs{a} \to 1$ for each
$n$. Proceeding inductively, we can use these properties of
$(f_n)$ to find increasing sequences of integers
$n_k \geq 1$ and numbers $0 < r_k < 1$ (starting with
$r_1 = \tfrac{1}{2}$, say) such that for each $k \geq 1$ one has
$\norm{f_{n_k}}_{H^2} < 2^{-k-1}$ and
\[
   \sup_{\abs{a} \leq r_k} \gamma(f_{n_k},a) < 2^{-k-1},  \qquad
   \sup_{\abs{a} \geq r_{k+1}} \gamma(f_{n_k},a) < 2^{-k-1}.
\]
For every $a \in \D$ we then have $\gamma(f_{n_k},a) < 2^{-k-1}$
for all except possibly one index $k$, for which
$\gamma(f_{n_k},a) \leq 1$. Hence
$\sum_k \gamma(f_{n_k},a) < 1 + \tfrac{1}{2} = \tfrac{3}{2}$.

Given a sequence $\lambda = (\lambda_k) \in c_0$, define
\[
   S\lambda = \sum_{k=1}^\infty \lambda_k f_{n_k}.
\]
The exponential decay of the $H^2$ norms of the functions $f_{n_k}$
ensures that the series converges in $H^2$ and hence pointwise. In
particular, from the fact that
$\abs{f_{n_k}(0)} \leq \norm{f_{n_k}}_{H^2} < 2^{-k-1}$
we get that $\abs{S\lambda(0)} \leq \tfrac{1}{2}\norm{\lambda}_\infty$.
In addition, for $a \in \D$,
\[
   \gamma(S\lambda,a)
   \leq \sum_{k=1}^\infty \abs{\lambda_k}\gamma(f_{n_k},a)
   \leq \frac{3}{2} \norm{\lambda}_\infty.
\]
Hence $\norm{S\lambda}_* \leq 2\norm{\lambda}_\infty$.
To check that $S\lambda \in \VMOA$, we let $\varepsilon > 0$ and
choose an integer $K$ such that $\abs{\lambda_k} \leq \varepsilon$
for $k > K$. Then, by estimating as above we have
\[
   \gamma(S\lambda,a)
   \leq \norm{\lambda}_\infty \sum_{k=1}^K \gamma(f_{n_k},a) +
        \frac{3}{2}\varepsilon.
\]
Since  $\gamma(f_{n_k},a) \to 0$ as $\abs{a} \to 1$ for each $k$, and
$\varepsilon > 0$ was arbitrary, this implies that
$S\lambda \in \VMOA$. Thus we have proved that
$S$ is a bounded linear operator from $c_0$ into $\VMOA$.

It remains to show that $S$ is bounded below. Given
$\lambda = (\lambda_k) \in c_0$, we first choose an index $K$ for
which $\abs{\lambda_K} = \norm{\lambda}_\infty$. Then we pick
a point $a \in \D$ such that
$\gamma(f_{n_K},a) \geq \tfrac{3}{4}$; this is possible
since $\norm{f_{n_K}}_* = 1$ and
$\abs{f_{n_K}(0)} < \tfrac{1}{4}$. Note that for $k \neq K$
we necessarily have $\gamma(f_{n_k},a) < 2^{-k-1}$.
Therefore, by employing the triangle inequality we get that
\[ \begin{split}
   \norm{S\lambda}_*
   &\geq \gamma(S\lambda,a)
   \geq \abs{\lambda_K} \gamma(f_{n_K},a) -
        \sum_{k\neq K} \abs{\lambda_k}\gamma(f_{n_k},a)  \\
   &\geq \frac{3}{4}\norm{\lambda}_\infty -
        \frac{1}{2}\norm{\lambda}_\infty
   = \frac{1}{4}\norm{\lambda}_\infty.
\end{split} \]
\end{proof}

\begin{remark} \label{re:Sub}
Let $X$ be a closed subspace of $\VMOA$.
As a consequence of Proposition~\ref{prop:c0} one has the
following dichotomy (see \cite{Leibov}):\ either $X$
contains an isomorphic copy of $c_0$ or the natural embedding
of $X$ into $H^2$ is an isomorphism. An analogous result in the setting
of martingale $\VMO$ spaces has been proved in \cite{MullerSchechtman}.
We point out here that the subspace structure of $\BMOA$ is very
complicated; see~\cite{Muller}.
\end{remark}

As noted at the beginning of the present section, the following
proposition completes the proof of Theorem~\ref{thm:Main1}.

\begin{prop} \label{prop:Comp}
Let $\varphi\colon \D \to \D$ be an analytic map and suppose that
condition \eqref{L} fails. Then the composition operator
$C_\varphi\colon \BMOA \to \BMOA$ fixes a copy of $c_0$ and therefore
it is not weakly compact.
\end{prop}

\begin{proof}
Since \eqref{L} fails to hold, we can find points $a_n \in \D$ such
that $\abs{\varphi(a_n)} \to 1$ and
\[
    \bignorm{\sigma_{\varphi(a_n)}\circ\varphi\circ\sigma_{a_n}}_{H^2}
    \geq c
\]
for some $c > 0$. Put $f_n = \sigma_{\varphi(a_n)}-\varphi(a_n)$. Then
$f_n(0) = 0$ and, for each $a \in \D$,
\[
    \bignorm{f_n\circ\sigma_a - f_n(a)}_{H^2}
    = \bignorm{\sigma_{\varphi(a_n)}\circ\sigma_a -
            \sigma_{\varphi(a_n)}(a)}_{H^2}
    = \sqrt{1-\abs{\sigma_{\varphi(a_n)}(a)}^2}.
\]
The last equality can be seen by using the fact that
$\sigma_{\varphi(a_n)}\circ\sigma_a$ is an inner function.
Now it follows easily that $f_n \in \VMOA$ and
$\norm{f_n}_* = 1$ for each $n$. By taking $a=0$ we 
obtain that $\norm{f_n}_{H^2} \to 0$ as $n \to \infty$.
Moreover,
\[
   \norm{C_\varphi f_n}_*
   \geq \bignorm{f_n\circ\varphi\circ\sigma_{a_n} -
              f_n(\varphi(a_n))}_{H^2}
   = \bignorm{\sigma_{\varphi(a_n)}\circ\varphi\circ\sigma_{a_n}}_{H^2}
   \geq c.
\]
According to Proposition~\ref{prop:c0} 
there is a subsequence $(f_{n_k})$ which is equivalent to the
natural basis of $c_0$. In particular, $(C_\varphi f_{n_k})$
is a weak-null sequence in $\BMOA$. By applying the
Bessaga-Pe{\l}czy\'{n}ski selection principle (see e.g.\
\cite[1.3.10]{AK}) to $(C_\varphi f_{n_k})$ we can pass to a
further subsequence, still denoted $(f_{n_k})$, such that
$(C_\varphi f_{n_k})$ is a semi-normalized basic sequence in $\BMOA$.
It follows that there are constants $A, B > 0$ so that 
\[
   A \cdot \norm{\lambda}_\infty
   \leq \biggnorm{ \sum_k \lambda_k C_\varphi f_{n_k} }_*
   \leq \norm{C_\varphi} \cdot \biggnorm{ \sum_k \lambda_k f_{n_k} }_*
   \leq B \cdot \norm{C_\varphi} \norm{\lambda}_\infty
\]
holds for any sequence $\lambda = (\lambda_k) \in c_0$. (To find $A$
just apply the biorthogonal basis functionals to
$\sum_k \lambda_k C_\varphi f_{n_k}$.) These estimates state that
the restriction of $C_\varphi$  to the closed subspace of $\BMOA$
spanned by the sequence $(f_{n_k})$ is an isomorphism
on a linearly isomorphic copy of $c_0$, and we are done.
\end{proof}

\begin{remark}
(1) Theorem \ref{thm:Main1} and its condition \eqref{L} 
also characterize the compactness, as well as the weak compactness,
of $C_\varphi$ on the space $\BMO$ identified with the space of
\emph{harmonic} functions $\D \to \C$ whose boundary values have bounded
mean oscillation. Indeed, it is known that a composition operator is
compact on $\BMOA$ if and only if it is compact on $\BMO$ (see e.g.\
\cite[Thm~3.5]{BCM}). Hence it remains to observe that if $C_\varphi$ is
weakly compact on $\BMO$, then it is weakly compact on the subspace
$\BMOA$ as well so that \eqref{L} holds.

(2) Theorem~\ref{thm:Main1} allows one to complete some
characterizations in \cite{Laitila:wcompact,Laitila:vector} as
follows:\ if $X$ is an infinite-dimensional complex reflexive Banach
space, then $C_\varphi$ is weakly compact on certain $X$-valued
versions of $\BMOA$ precisely when $C_\varphi$ is compact on $\BMOA$.
We refer to \cite{Laitila:wcompact,Laitila:vector} for a description
of this setting.
\end{remark}

\section{A condition on mean oscillation}
\label{sec:MeanOsc}

In this section our aim is to examine the function-theoretic meaning
of condition~\eqref{L} by revisiting the point of view that we already
touched upon in Section~\ref{sec:Cpt}. That is, \eqref{L} can be
thought of as a kind of pseudo-hyperbolic vanishing mean oscillation
condition for the boundary values of $\varphi$ over certain arcs
in $\T$; see Proposition~\ref{prop:Average} below.

To begin with we introduce some notation. When
$\varphi\colon \D \to \D$ is an analytic map and $I$ is an arc of $\T$,
denote 
\[
   \varphi_I = \frac{1}{\abs{I}} \int_I \varphi
   = \frac{1}{\abs{I}} \int_I \varphi(\zeta)\,\abs{d\zeta}
\]
for the integral average of $\varphi$ over $I$. Here and elsewhere
in this section all integrals over subsets of $\T$ are
calculated with respect to the normalized Lebesgue arc-length measure.
Also recall from Section~\ref{sec:Cpt} that
$I(re^{i\theta}) = \{ e^{it} : \abs{t-\theta} \leq \pi(1-r) \}$
for $re^{i\theta} \in \D$.

\begin{prop} \label{prop:Average}
For any analytic map $\varphi\colon \D \to \D$ condition~\eqref{L}
is equivalent to the following:
\begin{equation} \tag{A} \label{DACond}
   \frac{1}{\abs{I}^2} \int_I\int_I
   \rho\bigl( \varphi(\zeta),\varphi(\xi) \bigr)^2
   \,\abs{d\zeta}\abs{d\xi} \to 0
   \quad\text{as $\abs{\varphi_I} \to 1$},
\end{equation}
where  $I \subset \T$ are arcs.
\end{prop}

In the proof of this proposition we will make use of the following
easy estimate for the Poisson kernel, whose verification we leave
to the reader:\ for every $a \in \D$,
\begin{equation} \label{eq:PEst}
   \frac{1}{4\abs{I(a)}} \leq P_a(\zeta) \leq \frac{2}{\abs{I(a)}},
   \qquad \zeta \in I(a).
\end{equation}
We next record a simple auxiliary result, which isolates a crucial
step in proving Proposition~\ref{prop:Average}.

\begin{lemma} \label{le:Approach}
For $a \in \D$ we have $\abs{\varphi(a)} \to 1$ if and only if
$\abs{\varphi_{I(a)}} \to 1$.
\end{lemma}

\begin{proof}
The left-to-right implication is easy to prove. In fact, assuming
that $\varphi(a) \geq 0$ (as we may, after applying a rotation),
we get by using \eqref{eq:PEst} that
\[ \begin{split}
   1-\abs{\varphi(a)}
   =    \int_\T (1-\Re\varphi)P_a
   \geq \frac{1}{4\abs{I(a)}} \int_{I(a)} (1-\Re\varphi)
   \geq \frac{1}{4} \bigl(1-\abs{\varphi_{I(a)}}\bigr).
\end{split} \]
This clearly shows that $\abs{\varphi(a)} \to 1$ implies
$\abs{\varphi_{I(a)}} \to 1$.

For the reverse implication, we may assume that
$\varphi_{I(a)} \geq 1-\delta$ for some $0 < \delta < \tfrac{1}{2}$. Let
$E = \{ \zeta \in I(a) : \Re\varphi(\zeta) \geq 1-2\delta\}$.
Since $\Re\varphi \leq 1$, we must have
$\abs{E} \geq \tfrac{1}{2}\abs{I(a)}$. Consider the positive harmonic
function $u = \log(2/\abs{1-\varphi})$. It is geometrically obvious
that $\abs{1-\varphi} \leq c\sqrt{\delta}$ on $E$ for some constant
$c > 0$. Hence
\[
   u(a) \geq \int_\T uP_a
   \geq \biggl( \log\frac{2}{c\sqrt{\delta}} \biggr) \int_E P_a
   \geq \frac{1}{8} \biggl( \log\frac{2}{c\sqrt{\delta}} \biggr).
\]
Since $\abs{1-\varphi(a)} = 2e^{-u(a)}$, we deduce from this estimate
that $1-\abs{\varphi(a)} \leq \abs{1-\varphi(a)} \to 0$ as
$\delta \to 0$.
\end{proof}

\begin{proof}[Proof of Proposition~\ref{prop:Average}]
We start by proving the necessity of~\eqref{DACond}.  By the
preceding lemma $\abs{\varphi_I} \to 1$ implies that
$\abs{\varphi(a_I)} \to 1$. Hence \eqref{eq:POsc} and the left-hand
side of \eqref{eq:PEst} yield 
\begin{equation} \tag{A'}\label{ACond}
   \frac{1}{\abs{I}}
   \int_\T \rho\bigl(\varphi(\zeta),\varphi(a_I)\bigr)^2 \,\abs{d\zeta}
   \to 0 \quad \text{as $\abs{\varphi_I} \to 1$,}
\end{equation}
where $I \subset \T$ is an arc and $a_I \in \D$ is the unique point for
which $I = I(a_I)$. Then \eqref{DACond} is obtained from \eqref{ACond}
by a simple application of the triangle inequality
$\rho(\varphi(\zeta),\varphi(\xi)) \leq
\rho(\varphi(\zeta),\varphi(a_I)) + \rho(\varphi(\xi),\varphi(a_I))$.

To prove the sufficiency of \eqref{DACond} we will show that
\begin{equation} \label{eq:SuffAim}
   J(a) =
   \int_\T \int_\T \rho\bigl( \varphi(\zeta),\varphi(\xi) \bigr)^2
   P_a(\zeta) P_a(\xi)\,\abs{d\zeta}\abs{d\xi} \to 0
   \quad\text{as $\abs{\varphi(a)} \to 1$.}
\end{equation}
In view of \eqref{eq:POsc} this actually implies \eqref{L}, because
the function $w \mapsto \rho(z,w)^2$ is subharmonic in $\D$ and
therefore $\int_\T \rho\bigl(z,\varphi(\xi))^2 P_a(\xi)\,\abs{d\xi}
\geq \rho(z,\varphi(a))^2$ for every $z \in \overline{\D}$.

Let $\varepsilon > 0$. For each $a \in \D$ we can choose a point $a'$
on the line segment between $0$ and $a$ such that
$\int_{I(a')} P_a \geq 1-\varepsilon$ and
$1-\abs{a'} \leq c_\varepsilon (1-\abs{a})$ for some constant
$c_\varepsilon > 0$.  For real $a$ close to $1$ this can be seen
by integrating the estimate
$P_a(e^{it}) \geq (1-a^2)/[(1-a)^2+t^2]$ over an
interval $\abs{t} \leq c(1-a)$ and letting $c \to \infty$.
Now $\int_{\T\setminus I(a')} P_a \leq \varepsilon$, and since
$\rho \leq 1$, we can estimate
\[ \begin{split}
   J(a)
   &\leq 2\varepsilon + \int_{I(a')}\int_{I(a')}
      \rho\bigl( \varphi(\zeta),\varphi(\xi) \bigr)^2
      P_a(\zeta) P_a(\xi)\,\abs{d\zeta}\abs{d\xi}  \\
   &\leq 2\varepsilon +
      \frac{4c_\varepsilon^2}{\abs{I(a')}^2} \int_{I(a')}\int_{I(a')}
      \rho\bigl( \varphi(\zeta),\varphi(\xi) \bigr)^2
      \,\abs{d\zeta}\abs{d\xi}
\end{split} \]
by using the right-hand side of \eqref{eq:PEst} in the last step.
According to the Schwarz-Pick inequality we have
$\rho(\varphi(a),\varphi(a')) \leq \rho(a,a') \leq c_\varepsilon'$
for some $c_\varepsilon' < 1$ due to the fact that
$1-\abs{a'} \leq c_\varepsilon(1-\abs{a})$. Thus
$\abs{\varphi(a)} \to 1$ implies that $\abs{\varphi(a')} \to 1$, which,
in turn, yields $\abs{\varphi_{I(a')}} \to 1$ by
Lemma~\ref{le:Approach}.  By applying \eqref{DACond} to the arcs
$I(a')$ we then deduce from the above estimate that
$\limsup J(a) \leq 2\varepsilon$ as $\abs{\varphi(a)} \to 1$. Since
$\varepsilon > 0$ was arbitrary, this proves~\eqref{eq:SuffAim}.
\end{proof}

We summarize the principal function-theoretic compactness criteria for
$C_\varphi$ on $\BMOA$ in the following theorem. Criteria of a different
nature are given in \cite{BCM} and \cite{WirthsXiao}.

\begin{theorem}
Compactness and weak compactness of $C_\varphi\colon \BMOA \to \BMOA$
are equivalent to each of the conditions \eqref{Sm1}, \eqref{L},
\eqref{W1}, \eqref{W2}, \eqref{aver} and \eqref{ACond}.
\end{theorem}

\section{Results for VMOA}
\label{sec:VMOA}

In this section we discuss the case where $\varphi \in \VMOA$. Here 
simplified compactness criteria are available and new phenomena occur.
Recall first that if $\varphi \in \VMOA$ then $C_\varphi$ takes $\VMOA$
into itself and $C_\varphi\colon \BMOA \to \BMOA$ can be identified with
the biadjoint of its restriction to $\VMOA$; see
\cite[p.~939]{CimaMatheson}.

Let $\tau$ denote the hyperbolic metric in the unit disc, that is,
\[
   \tau(z,w) = \frac{1}{2} \log \frac{1+\rho(z,w)}{1-\rho(z,w)},
\]
where $\rho(z,w)$ is the pseudo-hyperbolic distance
between $z$ and $w$  (see e.g.\ \cite[Sec.~4.3]{Zhu}).
Contrary to the pseudo-hyperbolic metric, $\tau$ is unbounded
in $\D$ and it is appropriate to define $\tau (z,w)=\infty$ if
$z$ and $w$ are distinct points (at least) one of which lies on
the boundary.

We collect the main results in the case of $\VMOA$ as follows. 

\begin{theorem} \label{thm:vmoa}
Let $\varphi\colon \D \to \D$ be an analytic map such that
$\varphi \in \VMOA$. Then the following conditions are equivalent:
\begin{enumerate}
\renewcommand{\theenumi}{\roman{enumi}}
\item\label{it:CptVMOA}
$C_\varphi\colon \VMOA \to \VMOA$ is compact.
\item\label{it:WCptVMOA}
$C_\varphi\colon \VMOA \to \VMOA$ is weakly compact.
\item\label{it:JussiVMOA}
$\displaystyle
   \lim_{\abs{a} \to 1}
   \bignorm{\sigma_{\varphi(a)}\circ\varphi\circ\sigma_a}_{H^2} = 0$.
\item\label{it:phL}
$\displaystyle
   \lim_{\abs{a} \to 1}
   \int_\T \rho\bigl( \varphi(\sigma_a(\zeta)),\varphi(a) \bigr)^2
   \,\abs{d\zeta} = 0$.
\item\label{it:aver0}
$\displaystyle
   \lim_{\abs{I}\to 0} \frac{1}{\abs{I}^2} \int_I\int_I
   \rho\bigl( \varphi(\zeta),\varphi(\xi) \bigr)^2
   \,\abs{d\zeta}\abs{d\xi} = 0$,
   where $I \subset \T$ are arcs.
\end{enumerate}
Further, \eqref{it:phL} and \eqref{it:aver0} are equivalent to the
following conditions involving the hyperbolic metric:
\begin{enumerate}
\renewcommand{\theenumi}{\roman{enumi}'}
\setcounter{enumi}{3}
\item\label{it:MTj}
$\displaystyle
   \lim_{\abs{a} \to 1}
   \int_\T \tau\bigl( \varphi(\sigma_a(\zeta)),\varphi(a) \bigr)
   \,\abs{d\zeta} = 0$.
\item\label{it:avertau}
$\displaystyle
   \lim_{\abs{I}\to 0} \frac{1}{\abs{I}^2} \int_I\int_I
   \tau\bigl( \varphi(\zeta),\varphi(\xi) \bigr)
   \,\abs{d\zeta}\abs{d\xi} = 0$,
   where $I \subset \T$ are arcs.
\end{enumerate}
\end{theorem}

The main novelty of Theorem~\ref{thm:vmoa}, as compared to 
Theorem~\ref{thm:Main1}, lies in conditions \eqref{it:MTj}
and \eqref{it:avertau}, which relate to vanishing mean oscillation
with respect to the genuine hyperbolic metric. This also ties to
earlier research on composition operators from the Bloch space to
$\VMOA$. Before embarking on the proof of Theorem~\ref{thm:vmoa} we
discuss the interpretation of \eqref{it:MTj} from the literature and
draw some consequences.

First note that if the integral
$\int_\T \tau(\varphi(\sigma_a(\zeta)),\varphi(a))\,\abs{d\zeta}$
is finite for some $a \in \D$, then $\abs{\varphi} < 1$ a.e.\
on $\T$. Moreover, the integral stays bounded as $a$
varies on a compact subset of $\D$. Hence \eqref{it:MTj} implies 
\begin{equation}\label{CRU}
   \sup_{a \in \D} \int_\T
   \tau\bigl(\varphi(\sigma_a(\zeta)),\varphi(a)\bigr) \,\abs{d\zeta}
   < \infty,
\end{equation}
saying that $\varphi$ belongs to the hyperbolic $\BMOA$ class
introduced by Yamashita~\cite{Y}. Actually the fact that
\eqref{it:phL} implies the finiteness of the integral in \eqref{CRU}
for some $a \in \D$ is already non-trivial.

Recall that the Bloch space $\Bloch$ consists of the analytic functions
$f\colon \D \to \C$ for which
$\sup_{z\in\D} \abs{f'(z)}(1-\abs{z}^2) < \infty$. Then $\Bloch$
becomes a Banach space equipped with the norm
$\abs{f(0)} + \sup_{z\in\D} \abs{f'(z)}(1-\abs{z}^2)$. 
Composition operators $C_\varphi$ acting from $\Bloch$ into $\VMOA$ or
$\BMOA$ have been studied in e.g.\ \cite{Tjani}, \cite{CRU},
\cite{SZhao}, \cite{MakhTjani}, \cite{Xiao} and \cite{LMT}. As observed
by Makhmutov and Tjani \cite{MakhTjani}, it follows from the results of
Choe, Ramey and Ullrich \cite{CRU} combined with \cite{Y} that
$C_\varphi$ is bounded from $\Bloch$ into $\BMOA$ if and only if
\eqref{CRU} holds. In addition, it was proved in
\cite[Thm~6.1]{MakhTjani} that $C_\varphi$ is compact from $\Bloch$
into $\VMOA$ if and only if \eqref{it:MTj} holds. Therefore
Theorem~\ref{thm:vmoa} has the following surprising consequence.

\begin{cor} \label{cor:blochtovmoa}
Let $\varphi\colon \D \to \D$ be an analytic map with
$\varphi \in \VMOA$. Then $C_\varphi$ is compact $\VMOA \to \VMOA$
if and only if it is compact $\Bloch \to \VMOA$.
\end{cor}

This result was known earlier in the special case of boundedly valent
symbols $\varphi$ whose image $\varphi(\D)$ is contained in a polygon
inscribed in $\overline{\D}$; see \cite[Thm~5.3]{MakhTjani}. Of
course, in Corollary~\ref{cor:blochtovmoa} the implication from right
to left follows from the fact that $\VMOA$ is continuously embedded
in $\Bloch$. Furthermore, it is relevant to note that $C_\varphi$ is
bounded $\Bloch \to \VMOA$ if and only if it is compact
$\Bloch \to \VMOA$; see \cite[Thm~1.6]{SZhao}.

Towards the proof of Theorem~\ref{thm:vmoa} we make some
preliminary remarks. It was already observed by the first author
\cite[Thm~4.3]{Laitila:weighted} that condition \eqref{it:JussiVMOA}
alone characterizes the compactness of $C_\varphi\colon \VMOA \to \VMOA$.
At first sight \eqref{it:JussiVMOA} might seem stronger than 
\eqref{L} because $\abs{\varphi(a)} \to 1$ always implies
$\abs{a} \to 1$ by the Schwarz lemma. 
For the reader's convenience we include a direct function-theoretic
argument proving the equivalence of these two conditions
for symbols $\varphi \in \VMOA$.

\begin{lemma}\label{lm:jussi}
Let $\varphi\colon \D \to \D$ be an analytic map. Then condition
\eqref{it:JussiVMOA} of Theorem~\ref{thm:vmoa} holds if and only if
$\varphi \in \VMOA$ and \eqref{L} holds.
\end{lemma}

\begin{proof}
Let $\varphi_a = \sigma_{\varphi(a)}\circ\varphi\circ\sigma_a$.
By the self-inverse property of $\sigma_{\varphi(a)}$ we may write
$\varphi\circ\sigma_a = \sigma_{\varphi(a)}\circ\varphi_a$, from which
it follows that
\begin{equation} \label{eq:phi_a}
   \abs{(\varphi\circ\sigma_a)(z) - \varphi(a)} =
   \frac{1-\abs{\varphi(a)}^2}{\abs{1-\overline{\varphi(a)}\varphi_a(z)}}
   \abs{\varphi_a(z)}.
\end{equation}
This yields
$\norm{\varphi\circ\sigma_a - \varphi(a)}_{H^2} \leq
2\norm{\varphi_a}_{H^2}$. Hence \eqref{it:JussiVMOA} implies
that $\varphi \in \VMOA$.

Conversely note that if \eqref{L} holds but
\eqref{it:JussiVMOA} fails, then there exists a sequence
$(a_n)$ such that $\abs{a_n} \to 1$ while
$\abs{\varphi(a_n)} \leq r < 1$ and
$\norm{\varphi_{a_n}}_{H^2} \geq c > 0$ for all $n$. Then
\eqref{eq:phi_a} implies that
$\norm{\varphi\circ\sigma_{a_n} - \varphi(a_n)}_{H^2} \geq
(1-r)\norm{\varphi_{a_n}}_{H^2} \geq (1-r)c$, whence
$\varphi \notin \VMOA$. This proves the lemma.
\end{proof}

\begin{proof}[Proof of Theorem \ref{thm:vmoa}]
Recall that the operator $C_\varphi\colon \BMOA \to \BMOA$ is the
biadjoint of the restriction $C_\varphi\colon \VMOA \to \VMOA$,
since here $\varphi \in \VMOA$. Hence, according to
Theorem~\ref{thm:Main1}, conditions \eqref{it:CptVMOA} and
\eqref{it:WCptVMOA} are both equivalent to \eqref{L}. On the other hand,
in this case \eqref{L} and \eqref{it:JussiVMOA} are equivalent by
Lemma~\ref{lm:jussi}. We refer to Remark~\ref{rm:final} below
for an approach to the equivalences between conditions
\eqref{it:CptVMOA}--\eqref{it:JussiVMOA} which does not depend on
Section~\ref{sec:Cpt}.

Conditions \eqref{it:JussiVMOA} and \eqref{it:phL} are restatements of
each other according to \eqref{eq:POsc}. Furthermore, the equivalence
of \eqref{it:JussiVMOA} and \eqref{it:aver0} is proved in the same way
as Proposition~\ref{prop:Average}; instead of invoking
Lemma~\ref{le:Approach} we just observe that for points $a \in \D$
one has $\abs{a} \to 1$ if and only if $\abs{I(a)} \to 0$.

Since $\tau \geq c \rho^2$ for a suitable $c > 0$, it is obvious that
\eqref{it:avertau} implies \eqref{it:aver0}. Moreover,
\eqref{it:avertau} can be deduced from \eqref{it:MTj} by making a
change of variable, using the lower estimate from \eqref{eq:PEst} for
the Poisson kernel and applying the triangle inequality as in the
first part of the proof of Proposition~\ref{prop:Average}. 
The crucial remaining step in the
proof of Theorem~\ref{thm:vmoa} consists of verifying the
implication that the pseudo-hyperbolic condition \eqref{it:phL}
implies the hyperbolic condition \eqref{it:MTj}. We isolate this
more technical result below, which then completes the proof of the
theorem.
\end{proof}

\begin{prop} \label{prop:blochtobmoa}
Let $\varphi\colon \D \to \D$ be an analytic map.
Then condition \eqref{it:phL} implies condition \eqref{it:MTj}
in Theorem~\ref{thm:vmoa}.
\end{prop}

The argument will employ ideas of Wik \cite{Wik} related to his
elementary approach to the John-Nirenberg inequality for $\BMO$
functions. In particular, we will require the following one-dimensional
special case of \cite[Lemma~1]{Wik}:

\begin{lemma} \label{le:Wik}
Suppose that $0 < \lambda < 1$ and $E \subset [0,1]$ is any measurable
set having Lebesgue measure $\abs{E} \leq \lambda$. Then there is a
sequence $Q_1, Q_2, \ldots$ of closed dyadic intervals of $[0,1]$, 
having pairwise disjoint interiors, such that
$\tfrac{1}{2} \lambda \abs{Q_k} \leq \abs{Q_k \cap E}
\leq \lambda \abs{Q_k}$ for $k \geq 1$ and
$\bigabs{E \setminus \bigcup_k Q_k} = 0$.
\end{lemma}

\begin{proof}[Proof of Proposition \ref{prop:blochtobmoa}]
Assuming that condition \eqref{it:phL} (and equivalently also
\eqref{it:aver0}) holds, we split the proof into two steps. As the
first step we show:

\begin{claim} \label{claim:first}
$\displaystyle \lim_{\abs{a} \to 1} \frac{1}{\abs{I(a)}}
\int_{I(a)} \tau \bigl(\varphi(\zeta),\varphi(a)\bigr)
\,\abs{d\zeta} = 0$.
\end{claim}

To begin recall from Section~\ref{sec:Cpt} that condition
\eqref{it:phL} implies that $\abs{\varphi} < 1$ a.e.\ on $\T$ (this
fact can alternatively be deduced by observing that \eqref{it:CptVMOA}
implies the compactness of $C_\varphi$ on $H^2$ by \cite[Thm~4.1]{BCM}).
Towards the proof of Claim~\ref{claim:first} we first deduce from
\eqref{it:phL} by a change of variable and \eqref{eq:PEst} that
\begin{equation}\label{averMTj}
   \lim_{\abs{a} \to 1} \frac{1}{\abs{I(a)}} 
   \int_{I(a)} \rho \bigl(\varphi(\zeta),\varphi(a)\bigr)^2
   \,\abs{d\zeta} = 0,
\end{equation}
where $I(a) = \{e^{it}: \abs{t-\theta} \leq \pi(1-r)\}$ 
is the subarc of $\T$ associated to $a = re^{i\theta} \in \D$.
Hence we may pick $\delta > 0$ small enough so that
\begin{equation}\label{vika'}
   \frac{1}{\abs{I(a)}} \int_{I(a)}
   \rho \bigl(\varphi(\zeta),\varphi(a)\bigr)^2 \,\abs{d\zeta}
   < \frac{1}{4}
\end{equation}
whenever $a \in \D$ satisfies $\abs{a} > 1 - \delta$.

Let $\varepsilon \in (0,1/32)$. According to \eqref{it:aver0} 
we may decrease $\delta > 0$, if necessary, to ensure that for all
$a \in \D$ with $\abs{a} > 1 - \delta$ we also have
\begin{equation}\label{vika}
   \frac{1}{\abs{I(a)}^2} \int_{I(a)} \int_{I(a)}
   \rho \bigl(\varphi(\zeta),\varphi(\xi)\bigr)^2
   \,\abs{d\zeta}\abs{d\xi} < \varepsilon.
\end{equation}
Fix such a point $a$ and put
\[
   C_k = \bigl\{ \zeta \in I(a):
   \tau (\varphi(\zeta),\varphi(a)) \geq k \bigr\},
   \qquad k = 0,1,2,\ldots,
\]
whence $I(a) = C_0 \supset C_1 \supset C_2 \supset \cdots$.
Observe that if $\zeta \in C_1$, then the definition 
of the hyperbolic metric yields 
$\rho(\varphi(\zeta),\varphi(a)) \geq \beta$, where 
$\beta = \frac{e^2-1}{e^2+1} > 1/\sqrt{2}$. One gets from
\eqref{vika'} that 
\[
   \beta^2 \frac{\abs{C_1}}{\abs{I(a)}}
   \leq \frac{1}{\abs{I(a)}} \int_{I(a)}
        \rho \bigl(\varphi(\zeta),\varphi(a)\bigr)^2 \,\abs{d\zeta}
   < \frac{1}{4},
\]
whence $\abs{C_1} \leq \tfrac{1}{2} \abs{I(a)}$.

Let $k \geq 1$ be fixed.
Then we may apply Lemma~\ref{le:Wik} to the set
$C_k$ relative to the arc $I(a)$ with $\lambda = \tfrac{1}{2}$, which
gives a sequence $J_1, J_2, \ldots$ of subarcs of $I(a)$ with
disjoint interiors such that for each $\ell \geq 1$
\begin{equation}\label{arcs}
   \abs{C_k\cap J_\ell} \geq \tfrac{1}{4}\abs{J_\ell},
   \qquad
   \abs{C^c_k\cap J_\ell} \geq \tfrac{1}{2}\abs{J_\ell}
\end{equation}
and
\begin{equation} \label{density}
   \biggabs{C_k \setminus \bigcup_{\ell=1}^\infty J_\ell} = 0.
\end{equation}
Observe next that if $\zeta \in C_k^c$ and $\xi \in C_{k+1}$, then 
$\tau(\varphi(\zeta),\varphi(\xi)) \geq
\tau(\varphi(\xi),\varphi(a)) - \tau(\varphi(\zeta),\varphi(a)) \geq 1$,
so that $\rho(\varphi(\zeta),\varphi(\xi))^2 \geq \beta^2 > 0$. 
Consequently we get from \eqref{vika}, \eqref{arcs} and the assumption
on $a$ that 
\[ \begin{split}
   \varepsilon
   &> \frac{1}{\abs{J_\ell}^2} \int_{J_\ell} \int_{J_\ell} 
      \rho \bigl(\varphi(\zeta),\varphi(\xi)\bigr)^2
      \,\abs{d\zeta}\abs{d\xi}  \\
   &\geq \beta^2 \frac{\abs{C_k^c \cap J_\ell}}{\abs{J_\ell}}
      \cdot \frac{\abs{C_{k+1} \cap J_\ell}}{\abs{J_\ell}}
   \geq \frac{1}{4}
      \frac{\abs{C_{k+1} \cap J_\ell}}{\abs{J_\ell}}.
\end{split} \]
Thus $\abs{C_{k+1} \cap J_\ell} \leq 4\varepsilon \abs{J_\ell}$
for $\ell \geq 1$. We sum this inequality over $\ell$ and employ
\eqref{arcs} and \eqref{density} together with the essential
disjointness of the subarcs $J_\ell$ to obtain 
\begin{equation}\label{JNmeas}
   \abs{C_{k+1}}
   = \sum_{\ell=1}^\infty \abs{C_{k+1} \cap J_\ell}
   \leq 4\varepsilon \sum_{\ell=1}^\infty \abs{J_\ell}
   \leq 16\varepsilon \sum_{\ell=1}^\infty \abs{C_k \cap J_\ell}
   =    16\varepsilon \abs{C_k}.
\end{equation}
In particular, since $\varepsilon < 1/32$, we get by induction
that $\abs{C_k} \leq 2^{2-k}\abs{C_2}$ for $k \geq 2$.

Note that $k \leq \tau(\varphi(\zeta),\varphi(a)) < k + 1$ 
whenever $\zeta \in C_k \setminus C_{k+1}$ and $k \geq 0$.
Employing the short-hand notation $\{\tau < 2\}$ for the set
$\{\zeta \in I(a): \tau (\varphi(\zeta),\varphi(a)) < 2\} =
C_0 \setminus C_2$ we thus get that
\[ \begin{split}
   \int_{I(a)} \tau\bigl(\varphi(\zeta),\varphi(a)\bigr) \,\abs{d\zeta}
   &= \int_{\{\tau < 2\}}
      \tau\bigl(\varphi(\zeta),\varphi(a)\bigr) \,\abs{d\zeta} +
      \sum_{k=2}^\infty \int_{C_k \setminus C_{k+1}}
      \tau\bigl(\varphi(\zeta),\varphi(a)\bigr) \,\abs{d\zeta} \\
   &\leq \int_{\{\tau < 2\}}
      \tau\bigl(\varphi(\zeta),\varphi(a)\bigr) \,\abs{d\zeta} +
      \sum_{k=2}^\infty (k+1) \abs{C_k}.
\end{split} \] 
After division by $\abs{I(a)}$ the last term is less than
$\abs{C_2}\abs{I(a)}^{-1}\sum_{k=2}^\infty (k+1)2^{2-k}
\leq 128\varepsilon$, which tends to $0$ as $\varepsilon \to 0$.
On the other hand, in the set $\{\tau < 2\}$ we have 
$\tau(\varphi(\zeta),\varphi(a)) \leq
c\rho(\varphi(\zeta),\varphi(a))^2$ with a universal constant $c>0$,
so that also 
\[
   \lim_{\abs{a}\to 1} \frac{1}{\abs{I(a)}} \int_{\{\tau < 2\}}
   \tau \bigl(\varphi(\zeta),\varphi(a)\bigr) \,\abs{d\zeta} = 0
\]
in view of \eqref{averMTj}. This finishes the proof
of Claim~\ref{claim:first}.

\medskip

As the final step we show that the condition of Claim~\ref{claim:first}
implies the desired hyperbolic condition \eqref{it:MTj} of
Theorem~\ref{thm:vmoa}. The required argument is quite standard but 
more technical than the analogous fact for the
pseudo-hyperbolic distance $\rho$ in Section~\ref{sec:MeanOsc}
because the hyperbolic distance $\tau$ is unbounded. We omit some
computational details.

\begin{claim} \label{claim:second}
$\displaystyle
\int_\T \tau\bigl( \varphi(\sigma_a(\zeta)),\varphi(a) \bigr)
\,\abs{d\zeta} =
\int_\T \tau \bigl(\varphi(\zeta),\varphi(a)\bigr) P_a(\zeta)
\,\abs{d\zeta} \to 0$ as $\abs{a} \to 1$.
\end{claim}

For the proof we assume that $a \in \D$ satisfies
$2^{-N}\leq 1 - \abs{a} < 2^{1-N}$ 
for some $N \geq 1$, and then let $N \to \infty$ in our estimates. 
Define for $k = 1,\ldots,N$ the radii $r_k$, points $a_k\in \D$ and
arcs $I_k$ through  $1-r_k = 2^{N-k}(1-\abs{a})$,
$a_k=r_ka/\abs{a}$ and $I_k=I(a_k).$ Set also $a_0=0$
and $I_0=\T$. Then $a=a_N$ and
$I(a) = I_N\subset I_{N-1}\subset\ldots \subset I_0=\T.$
Moreover, $2^{-k}\leq \abs{I_k} < 2^{1-k}.$
Observe that if $1 \leq k<N$ and $\zeta \in I_k \setminus I_{k+1}$,
then elementary trigonometry yields
$\abs{\zeta-a} \geq \frac{1}{2} \abs{I_{k+1}} \geq 2^{-k-2}$.
Hence the Poisson kernel satisfies $P_a(\zeta) \lesssim 2^{2k-N}$
for all $\zeta \in I_k \setminus I_{k+1}$, where $\lesssim$ indicates
that the left-hand side is bounded above by a constant multiple
of the right-hand side, the constant being independent of $N$ and $k$.
Consequently we may estimate the second integral appearing in
Claim~\ref{claim:second} as follows: 
\[ \begin{split}
   &\int_\T \tau \bigl(\varphi(\zeta),\varphi(a)\bigr) P_a(\zeta)
    \,\abs{d\zeta} \\
   &\qquad\lesssim
    \sum_{k=0}^{N-1} 2^{2k-N} \int_{I_k\setminus I_{k+1}}
    \tau \bigl(\varphi(\zeta),\varphi(a)\bigr) \,\abs{d\zeta} +
    2^N \int_{I(a)} \tau \bigl(\varphi(\zeta),\varphi(a)\bigr)
        \,\abs{d\zeta} \\
   &\qquad\lesssim
    \sum_{k=0}^N  \frac{2^{k-N}}{\abs{I_k}}
    \int_{I_k} \tau \bigl(\varphi(\zeta),\varphi(a)\bigr)
    \,\abs{d\zeta} \\
   &\qquad\leq
    \sum_{k=0}^N \frac{2^{k-N}}{\abs{I_k}}
    \int_{I_k} \tau \bigl(\varphi(\zeta),\varphi(a_k)\bigr)
    \,\abs{d\zeta} +
    \sum_{k=0}^{N-1} 2^{k-N} \tau\bigl(\varphi(a_k),\varphi(a)\bigr) \\
   &\qquad\equiv A_N + B_N. 
\end{split} \]

It will suffice to verify that the condition of Claim~\ref{claim:first}
implies that the terms $A_N$ and $B_N$ both tend to zero as
$N \to \infty$. First of all (observe that now \eqref{CRU} holds), 
\[ \begin{split}
   A_N &\lesssim \Biggl( \sum_{k=0}^{[N/2]} 2^{k-N} +
      \sum_{k=[N/2]+1}^N 2^{k-N} \Biggr) \frac{1}{\abs{I_k}}
      \int_{I_k} \tau \bigl(\varphi(\zeta),\varphi(a_k)\bigr)
      \,\abs{d\zeta} \\
     &\lesssim N \cdot 2^{-N/2} + \sup_{k>[N/2]} \frac{1}{\abs{I_k}}
      \int_{I_k} \tau \bigl(\varphi(\zeta),\varphi(a_k)\bigr)
      \,\abs{d\zeta}
\end{split} \]
Above the first term tends to zero trivially, and the second term by
Claim~\ref{claim:first}, as $N \to \infty$.

In order to relate the term $B_N$ to the averages in
Claim~\ref{claim:first} we introduce the short-hand
$b_k = \abs{I_k}^{-1}\int_{I_k} \tau (\varphi(\zeta),\varphi(a_k))
\,\abs{d\zeta}$. Let $1 \leq k\leq N$. By averaging over the arc
$I_k$ we get from the triangle inequality for $\tau$ that
\[ \begin{split}
   \tau \bigl(\varphi(a_{k-1}),\varphi(a_k)\bigr)
   &\leq \frac{1}{\abs{I_k}} \int_{I_k}
         \tau \bigl(\varphi(\zeta),\varphi(a_{k-1})\bigr)
         \,\abs{d\zeta} +
         \frac{1}{\abs{I_k}} \int_{I_k}
         \tau \bigl(\varphi(\zeta),\varphi(a_k)\bigr) \,\abs{d\zeta} \\
   &\leq 2b_{k-1} + b_k,
\end{split} \]
since $\abs{I_{k-1}} \leq 2\abs{I_k}$. Because $a = a_N$, we deduce
that
\[
   \tau \bigl(\varphi(a_k),\varphi(a)\bigr)
   \lesssim \sum_{j=k}^N b_j
   \le (N-k+1) \max_{k \leq j \leq N} b_j.
\]
Put $E_k = \max_{k \leq j \leq N} b_j$, so that by combining the above
estimates one has
\[
   B_N \lesssim \sum_{k=0}^{N-1} (N - k+1)2^{k-N} E_k,
\]
where the $E_k$'s have a uniform upper bound (independent of $a$) and
$E_{[N/2]} \to 0$ as $N \to \infty$. By splitting the preceding sum as
before at the level $[N/2]$ we deduce that $B_N \to 0$ as
$N \to \infty$. This completes the proof of Claim~\ref{claim:second},
and hence of Proposition~\ref{prop:blochtobmoa}.
\end{proof}

\begin{remark}\label{rm:final}
(1) The equivalence of conditions
\eqref{it:CptVMOA}--\eqref{it:JussiVMOA} in Theorem~\ref{thm:vmoa} can
be proved without relying on the work of Section~\ref{sec:Cpt}.
One essentially argues as in the proof of Proposition \ref{prop:Comp} and
invokes Lemma~\ref{lm:jussi} together with the comments preceding it.
Instead of using the Bessaga-Pe{\l}czy\'{n}ski selection principle
one may just apply Proposition~\ref{prop:c0} twice, the second time
to the image sequence. We leave the
details to the interested reader.

(2) In \cite{MakhTjani} an analytic map $\varphi\colon \D \to \D$ is
said to belong to the hyperbolic class $\VMOA^h$ if $\varphi$ satisfies
\eqref{it:MTj}. Similarly, we may say that $\varphi$ belongs to the
pseudo-hyperbolic class $\VMOA^{ph}$ if \eqref{it:phL} holds.
Thus Proposition~\ref{prop:blochtobmoa} (and its converse) states that
$\varphi \in \VMOA^h$ if and only if $\varphi \in \VMOA^{ph}$, 
which is an interpretation independent of composition operators.

(3) 
In the formulation of conditions \eqref{it:MTj} and \eqref{it:phL} the
metrics $\rho$ and $\tau$ are raised to different powers. However, in
each condition the power is irrelevant. Namely, an inspection of the
proof of Proposition~\ref{prop:blochtobmoa} shows that one may replace
$\tau$ by any power $\tau^p$ with $p>0$ in \eqref{it:MTj}. This yields
the same conclusion for condition \eqref{it:avertau}, and the analogous
fact for \eqref{it:phL} and \eqref{it:aver0} is obvious.

(4) Proposition~\ref{prop:blochtobmoa} suggests the following problem, 
which we did not pursue any further:\ is there a version of the
proposition for composition operators $\Bloch \to \BMOA$?
We recall here that Xiao \cite{Xiao} (cf.\ also \cite{LMT}) showed
that $C_\varphi$ is compact $\Bloch \to \BMOA$
if and only if 
\[
   \lim_{r\to 1} \sup_{a \in \D} \int_{\{z: \abs{\varphi(z)} > r\}}
   \frac{\abs{\varphi'(z)}^2}{(1-\abs{\varphi(z)}^2)^2}
   (1-\abs{\sigma_a(z)}^2) \,dA(z) = 0,
\]
where $A$ is the planar Lebesgue measure.
\end{remark}

\bibliographystyle{amsalpha}

\end{document}